\documentclass[a4paper,10pt,leqno]{article}
\usepackage{amsmath,amsfonts,amssymb,amsthm}
\usepackage{graphicx,subfigure,booktabs,multirow}
\usepackage{mathrsfs,enumerate,color,bm}
\usepackage{algorithm,algpseudocode}
\usepackage{charter}
\usepackage{geometry}
\newtheorem{theorem}{Theorem}[section]

\newtheorem{proposition}[theorem]{Proposition}

\title{On Some Sufficient Conditions for Strong Ellipticity}
\author{
Weiyang Ding\thanks{Department of Applied Mathematics, The Hong Kong Polytechnic University, Kowloon, Hong Kong. Email: {\tt weiyang.ding@gmail.com}. This author¡¦s work was partially supported by the Hong Kong Research Grants Council (Grant No. C1007-15G).}
\and
Liqun Qi\thanks{Department of Applied Mathematics, The Hong Kong Polytechnic University, Kowloon, Hong Kong. Email: {\tt liqun.qi@polyu.edu.hk}. This author¡¦s work was partially supported by the Hong Kong Research Grants Council (Grant No. PolyU 501913, 15302114, 15300715, 15301716 and C1007-15G).}
\and
Hong Yan\thanks{Department of Electronic Engineering, City University of Hong Kong, Kowloon, Hong Kong. Email: {\tt h.yan@cityu.edu.hk}. This author¡¦s work was partially supported by the Hong Kong Research Grants Council (Grant No. C1007-15G).}
}
\date{\today}

\begin{document}

\maketitle

\begin{abstract}
  We establish several sufficient conditions for the strong ellipticity of any fourth-order elasticity tensor in this paper.
  The first presented sufficient condition is an extension of positive definite matrices, which states that the strong ellipticity holds if the unfolding matrix of this fourth-order elasticity tensor can be modified into a positive definite one by preserving the summations of some corresponding entries.
  An alternating projection algorithm is proposed to verify whether an elasticity tensor satisfies the first condition or not.
  Conditions for some special cases beyond the first sufficient condition are further investigated, which includes some important cases for the isotropic and some particular anisotropic linearly elastic materials.



\end{abstract}

\bigskip

\section{Introduction}

The strong ellipticity condition is known for a long time as one of the basic constitutive ingredients in the theory of elasticity, which guarantees the existence of solutions of basic boundary-value problems of elastostatics and thus ensures an elastic material to satisfy some mechanical properties, see eg. \cite{Aron83,ChiritaGhiba10,GourgiotisBigoni16,Gurtin73,KnowlesSternberg75,KnowlesSternberg76,Rosakis90} and works cited therein.
Thereby, it is essential to identify whether the strong ellipticity holds or not in subdomains of the domain of the strain energy function.

The strong ellipticity of isotropic materials is relatively easier to handle, thus the early related works focus on isotropic materials.
In particular cases, it is well-known that the strong ellipticity of an isotropic linear elastic material is equivalent to two simple inequalities about the Lam\'{e} moduli \cite{Gurtin73}.
Knowles and Sternberg \cite{KnowlesSternberg75,KnowlesSternberg76} established necessary and sufficient conditions for both ordinary and strong ellipticity of the equations governing finite plane equilibrium deformations of a compressible hyperelastic solid.
Simpson and Spector \cite{SimpsonSpector} extended their works to the spacial case using the representation theorem for copositive matrices.
Some reformulations were also given in Rosakis \cite{Rosakis90} and Wang and Aron \cite{WangAron96}.
One can refer to \cite{Dacorogna01} as a review of necessary and sufficient conditions for strong ellipticity for isotropic cases.

As to the strong ellipticity of anisotropic materials, Walton and Wilber \cite{WaltonWilber} provided sufficient conditions for strong ellipticity of a general class of anisotropic hyperelastic materials, which require the
first partial derivatives of the reduced-stored energy function to satisfy several simple inequalities and the second partial derivatives to satisfy a convexity condition.
Chiri\c t\u a, Danescu, and Ciarletta\cite{ChiritaDanescuCiarletta07} and Zubov and Rudev \cite{ZubovRudev16} proposed sufficient and necessary conditions for the strong ellipticity of certain classes of anisotropic linearly elastic materials.
Gourgiotis and Bigoni \cite{GourgiotisBigoni16} investigated the strong ellipticity of materials with extreme mechanical anisotropy.

Qi, Dai, and Han \cite{QiDaiHan09} introduced M-eigenvalues for ellipticity tensors and proved that the strong ellipticity holds if and only if all the M-eigenvalues of the ellipticity tensor is positive.
Wang, Qi, and Zhang \cite{WangQiZhang09} proposed a practical power method for computing the largest M-eigenvalue of any ellipticity tensor, which can also be employed to verify the strong ellipticity.
Very recently, Huang and Qi \cite{HuangQi17} generalized the M-eigenvalues of fourth-order ellipticity tensors and related algorithms to higher order cases.
Chang, Qi, and Zhou \cite{ChangQiZhou10} defined another type of ``eigenvalues'' for ellipticity tensors named as singular values, and the positivity of all the singular values of the ellipticity tensor is also a necessary and sufficient condition for the strong ellipticity.
Han, Dai, and Qi \cite{HanDaiQi09} linked the strong ellipticity condition to the rank-one positive definiteness of three second-order tensors, three fourth-order tensors, and a sixth-order tensor.

The present paper is built up as follows. We briefly introduce the strong ellipticity for linearly elastic materials and some related concepts and notations in Section \ref{sec_se}.
Then in Section \ref{sec_convex}, we establish a sufficient condition for strong ellipticity by identifying whether the intersection of two closed convex sets are nonempty, and propose an alternating projection verification algorithm.
In Section \ref{sec_special}, we focus on special cases beyond the first sufficient condition, which includes important cases for the isotropic and particular anisotropic linearly elastic materials.
Finally, we draw concluding remarks and propose open questions in Section \ref{sec_conclusion}.

\section{Strong ellipticity for linearly elastic materials}\label{sec_se}

For a linearly elastic material, the tensor of elastic moduli in a Cartesian coordinate system is a fourth-order three-dimensional tensor $\mathscr{A} = (a_{ijkl}) \in \mathbb{R}^{3 \times 3 \times 3 \times 3}$ which is invariant under the following permutations of indices
\begin{equation}\label{eq_sym}
  a_{ijkl} = a_{jikl} = a_{ijlk}.
\end{equation}
We use $\mathbb{E}$ to denote the set of all fourth-order three-dimensional tensors satisfying \eqref{eq_sym}.
The strong ellipticity condition (SE-condition) is stated by
\begin{equation}\label{eq_se}
  \mathscr{A} {\bf x}^2 {\bf y}^2 := \sum_{i,j,k,l=1}^3 a_{ijkl} x_i x_j y_k y_l > 0
\end{equation}
for any nonzero vectors ${\bf x}, {\bf y} \in \mathbb{R}^3$.
We also call a tensor satisfying the SE-condition to be M-positive definite (M-PD) \cite{QiLuo17}.
Similarly, a tensor $\mathscr{A} \in \mathbb{E}$ is said to be M-positive semidefinite (M-PSD) \cite{QiLuo17} if $\mathscr{A} {\bf x}^2 {\bf y}^2 \geq 0$ for any vectors ${\bf x}, {\bf y} \in \mathbb{R}^3$.

We can introduce another type of positive definiteness for fourth-order tensors, which is often involved in the investigations of linearly elastic materials. Let $\mathscr{A}$ be a fourth-order three-dimensional tensor and ${\bf Z}$ be a three-by-three matrix. Define
\begin{equation}
  \mathscr{A} {\bf Z}^2 := \sum_{i,j,k,l = 1}^3 a_{ijkl} z_{ik} z_{jl}.
\end{equation}
If $\mathscr{A} {\bf Z}^2 > 0$ ($\geq 0$) for any nonzero matrix ${\bf Z} \in \mathbb{R}^{3 \times 3}$, then $\mathscr{A}$ is said to be S-positive (semi)definite \cite{QiLuo17}.
Denote ${\bf A}_{kl} = (a_{ijkl}) \in \mathbb{R}^{3 \times 3}$ for each $k,l$ and ${\bf z}_k = (z_{ik}) \in \mathbb{R}^3$ for each $k$.
Then the tensor $\mathscr{A}$ and the matrix ${\bf Z}$ can be unfolded into a matrix ${\bf A} \in \mathbb{R}^{9 \times 9}$ and a vector ${\bf z} \in \mathbb{R}^9$ as
$$
{\bf A} =
\begin{bmatrix}
  {\bf A}_{11} & {\bf A}_{12} & {\bf A}_{13} \\
  {\bf A}_{21} & {\bf A}_{22} & {\bf A}_{23} \\
  {\bf A}_{31} & {\bf A}_{32} & {\bf A}_{33}
\end{bmatrix}
\quad {\rm and} \quad
{\bf z} =
\begin{bmatrix}
  {\bf z}_1 \\ {\bf z}_2 \\ {\bf z}_3
\end{bmatrix},
$$
respectively.
Moreover, it can be verified that $\mathscr{A} {\bf Z}^2 = {\bf z}^\top {\bf A} {\bf z}$.
Therefore the fourth-order tensor $\mathscr{A}$ is S-PD or S-PSD if and only if the matrix ${\bf A}$ is PD or PSD, respectively.

The S-positive definiteness is a sufficient condition for the M-positive definiteness, i.e., the SE-condition, which can be proved by observing $\mathscr{A} {\bf Z}^2 = \mathscr{A} {\bf x}^2 {\bf y}^2$ when ${\bf Z} = {\bf x} {\bf y}^\top$.
Nevertheless, the S-positive definiteness is not a necessary condition for the SE-condition.
A counter example is a tensor with
$$
a_{1111} = a_{2222} = a_{3333} = 2, \quad a_{1221} = a_{2121} = a_{2112} = a_{1212} = 1,
$$
and all other entries equal to zero.
Then the bi-quadratic form is
$$
\mathscr{A} {\bf x}^2 {\bf y}^2 = 2(x_1 y_1 + x_2 y_2)^2 + 2 x_3^2 y_3^2,
$$
thus $\mathscr{A}$ is apparently M-PSD, while the unfolding matrix
$$
\left[
\begin{array}{ccc|ccc|ccc}
  2 & 0 & 0 & 0 & 1 & 0 & 0 & 0 & 0 \\
  0 & 0 & 0 & 1 & 0 & 0 & 0 & 0 & 0 \\
  0 & 0 & 0 & 0 & 0 & 0 & 0 & 0 & 0 \\
  \hline
  0 & 1 & 0 & 0 & 0 & 0 & 0 & 0 & 0 \\
  1 & 0 & 0 & 0 & 2 & 0 & 0 & 0 & 0 \\
  0 & 0 & 0 & 0 & 0 & 0 & 0 & 0 & 0 \\
  \hline
  0 & 0 & 0 & 0 & 0 & 0 & 0 & 0 & 0 \\
  0 & 0 & 0 & 0 & 0 & 0 & 0 & 0 & 0 \\
  0 & 0 & 0 & 0 & 0 & 0 & 0 & 0 & 2
\end{array}
\right]
$$
is not PSD.
Therefore, we desire to explore weaker but still easily checkable sufficient conditions for the SE-condition.

\section{A sufficient condition with a verification algorithm}\label{sec_convex}

For any $\mathscr{A} \in \mathbb{E}$ and ${\bf y} \in \mathbb{R}^3$, we denote a three-by-three matrix $\mathscr{A} {\bf y}^2$ as
$$
(\mathscr{A} {\bf y}^2)_{ij} := \sum_{k,l=1}^n a_{ijkl} y_k y_l, \quad i,j = 1,2,3.
$$
Since $\mathscr{A}$ admits the symmetries in \eqref{eq_sym}, this matrix $\mathscr{A} {\bf y}^2$ is a symmetric matrix.
Noticing that $\mathscr{A} {\bf x}^2 {\bf y}^2 = {\bf x}^\top (\mathscr{A} {\bf y}^2) {\bf x}$ for any vectors ${\bf x}, {\bf y} \in \mathbb{R}^3$, we can prove the following necessary and sufficient condition for the M-positive (semi)definiteness.
\begin{proposition}
  Let $\mathscr{A} \in \mathbb{E}$.
  Then $\mathscr{A}$ is M-PD or M-PSD if and only if the matrix $\mathscr{A} {\bf y}^2$ is PD or PSD for each nonzero ${\bf y} \in \mathbb{R}^3$, respectively.
\end{proposition}

Generally speaking, the above necessary and sufficient condition is as hard as the SE-condition to check.
However, it motivates some checkable sufficient conditions.
Recall that every positive semidefinite matrix can be decomposed into the sum of rank-one positive semidefinite matrices and the minimal number of terms is exactly its rank \cite{HornJohnson13}.
Thus, we have the following sufficient condition for a tensor $\mathscr{A}$ to be M-PSD 
\begin{equation}\label{eq_rank1+}
  \mathscr{A} {\bf y}^2 = \sum_{s = 1}^r \alpha_s {\bf f}_s({\bf y}) {\bf f}_s({\bf y})^\top, \quad \alpha_s > 0,
\end{equation}
where each ${\bf f}_s({\bf y})$ is a homogeneous function of degree one, i.e., ${\bf f}_s({\bf y}) = {\bf U}_s {\bf y}$ for $s = 1,2,\dots,r$.
Any matrix $\mathscr{A} {\bf y}^2$ in the above form is PSD and thus $\mathscr{A}$ is M-PSD.
Furthermore, if $\mathscr{A}$ is M-PD then the number of terms in the summation should be no less than $3$, i.e., $r \geq 3$.
Denote the entries of each ${\bf U}_s$ as $u_{ij}^{(s)}$ ($i,j=1,2,3$).
Then \eqref{eq_rank1+} reads
\[
\begin{split}
\sum_{k,l=1}^n a_{ijkl} y_k y_l
&= \sum_{s=1}^r \alpha_s \Big( \sum_{k=1}^n u_{ik}^{(s)} y_k \Big) \Big( \sum_{l=1}^n u_{jl}^{(s)} y_l \Big) \\
&= \sum_{k,l=1}^n \Big( \sum_{s=1}^r \alpha_s u_{ik}^{(s)} u_{jl}^{(s)} \Big) y_k y_l.
\end{split}
\]
Therefore, given ${\bf U}_s$ ($s = 1,2,\dots,r$), the entries of $\mathscr{A}$ are uniquely determined by
\begin{equation}\label{eq_con}
  a_{ijkl} = \frac{1}{2} \sum_{s=1}^r \alpha_s \Big( u_{ik}^{(s)} u_{jl}^{(s)} + u_{jk}^{(s)} u_{il}^{(s)} \Big),
\end{equation}
which satisfies the symmetries $a_{ijkl} = a_{jikl} = a_{ijlk}$.

Next, we shall discuss when a tensor in $\mathbb{E}$ can be represented by \eqref{eq_con}.
Denote another fourth-order three-dimensional tensor $\mathscr{B}$ by
$$
b_{ijkl} = \sum_{s=1}^r \alpha_s u_{ik}^{(s)} u_{jl}^{(s)}.
$$
Note that $\mathscr{B}$ may not in the set $\mathbb{E}$, i.e., it is not required to obey \eqref{eq_sym}, but it still satisfies a weaker symmetry that $b_{ijkl} = b_{jilk}$.
It can be seen that its unfolding ${\bf B}$ is a PSD matrix from
$$
{\bf B} = \sum_{s=1}^r \alpha_s {\bf u}_s {\bf u}_s^\top,
$$
where ${\bf u}_s$ is the unfolding (or vectorization) of ${\bf U}_s$ ($s = 1,2,\dots,r$).
Hence $\mathscr{B}$ is S-PSD since all the coefficients $\alpha_s$ are positive.
Furthermore, comparing the entries of $\mathscr{A}$ and $\mathscr{B}$, we will find that
$$
a_{ijkl} = a_{jikl} = \frac{1}{2} (b_{ijkl} + b_{jikl}), \quad i,j,k,l = 1,2,3,
$$
and thus
$$
\mathscr{A} {\bf x}^2 {\bf y}^2 = \mathscr{B} {\bf x}^2 {\bf y}^2 = \mathscr{B} \big( {\bf x} {\bf y}^\top \big)^2.
$$
Therefore, $\mathscr{A}$ is M-PD or M-PSD when $\mathscr{B}$ is S-PD or S-PSD, respectively.

Given a fourth-order tensor $\mathscr{A} \in \mathbb{E}$, we denote
$$
\mathbb{T}_{\mathscr{A}} := \{ \mathscr{T}:\, t_{ijkl} = t_{jilk},\, t_{ijkl} + t_{jikl} = 2 a_{ijkl} \}.
$$
We also denote the set of all fourth-order S-PSD tensors as
$$
\mathbb{S} := \{ \mathscr{T}:\, t_{ijkl} = t_{jilk},\, \mathscr{T} \text{ is S-PSD} \}.
$$
Note that both $\mathbb{T}_{\mathscr{A}}$ and $\mathbb{S}$ are closed convex sets, where $\mathbb{T}_{\mathscr{A}}$ is a linear subspace of the whole space of all the fourth-order three-dimensional tensor with $t_{ijkl} = t_{jilk}$ and $\mathbb{S}$ is isomorphic with the nine-by-nine symmetric PSD matrix cone.
Furthermore, we have actually proved the following sufficient condition for a tensor to be M-PD or M-PSD.
\begin{theorem}
  Let $\mathscr{A} \in \mathbb{E}$.
  If $\mathbb{T}_{\mathscr{A}} \cap \mathbb{S} \neq \emptyset$, then $\mathscr{A}$ is M-PSD;
  If $\mathbb{T}_{\mathscr{A}} \cap (\mathbb{S} \setminus \partial\mathbb{S}) \neq \emptyset$, then $\mathscr{A}$ is M-PD.
\end{theorem}

A method called projections onto convex sets (POCS) \cite{BauschkeBorwein96,EscalanteRaydan11} is often employed to check whether the intersection of two closed convex sets is empty or not.
POCS is also known as the alternating projection algorithm.
Denote ${\cal P}_1$ and ${\cal P}_2$ as the projections onto $\mathbb{T}_{\mathscr{A}}$ and $\mathbb{S}$, respectively.
Then POCS is stated by
$$
\left\{
\begin{array}{l}
  \mathscr{B}^{(t+1)} = {\cal P}_2( \mathscr{A}^{(t)} ), \\
  \mathscr{A}^{(t+1)} = {\cal P}_1( \mathscr{B}^{(t+1)} ),
\end{array}
\right.
\quad t = 0,1,2,\dots.
$$
The algorithm can be described as the following iterative scheme:
\begin{equation}\label{alg_pocs}
\left\{
\begin{array}{l}
  \text{Eigendecomposition } {\bf A}^{(t)} = {\bf V}^{(t)} {\bf D}^{(t)} ({\bf V}^{(t)})^\top, \\
  {\bf B}^{(t+1)} = {\bf V}^{(t)} {\bf D}_+^{(t)} ({\bf V}^{(t)})^\top, \\
  a_{iikl}^{(t+1)} = a_{iikl} \text{ for } i,k,l = 1,2,3, \\
  a_{ijkk}^{(t+1)} = a_{ijkk} \text{ for } i,j,k = 1,2,3, \\
  a_{ijkl}^{(t+1)} = a_{ijkl} + \frac{1}{2} (b_{ijkl}^{(t+1)} - b_{jikl}^{(t+1)}) \text{ for } i \neq j,\ k \neq l,
\end{array}
\right.
\quad t = 0,1,2,\dots.
\end{equation}
where $\mathscr{A}^{(0)} = \mathscr{A}$, ${\bf A}^{(t)}$ and ${\bf B}^{(t)}$ are the unfolding matrices of $\mathscr{A}^{(t)}$ and $\mathscr{B}^{(t)}$ respectively, and ${\bf D}_+^{(t)} = {\rm diag}\big( \max(d_{ii}^{(t)}, 0) \big)$.
The convergence of the alternating projection method between two closed convex sets has been known for a long time \cite{CheneyGoldstein59}.

\begin{theorem}
  Let $\mathscr{A} \in \mathbb{E}$. If $\mathbb{T}_{\mathscr{A}} \cap \mathbb{S} \neq \emptyset$, then the sequences $\{ \mathscr{A}^{(t)} \}$ and $\{ \mathscr{B}^{(t)} \}$ produced by Algorithm \eqref{alg_pocs} both converge to a point $\mathscr{A}^\ast \in \mathbb{T}_{\mathscr{A}} \cap \mathbb{S}$.
\end{theorem}

Because the convergence of POCS requires the involved convex sets to be closed, Algorithm \eqref{alg_pocs} is only suitable for identifying the M-positive semidefiniteness.
If we want to check the M-positive definiteness, then some modifications are needed.
Denote $\mathscr{E} \in \mathbb{E}$ with $e_{iikk} = 1$ ($i,k=1,2,3$) and other entries being zero.
Note that $\mathscr{E} {\bf x}^2 {\bf y}^2 = ({\bf x}^\top {\bf x}) ({\bf y}^\top {\bf y})$.
Hence $\mathscr{E}$ is M-PD, which implies that $\mathscr{A}$ is M-PD if and only if $\mathscr{A} - \epsilon \mathscr{E}$ is M-PSD for some sufficiently small $\epsilon > 0$.
From such observation, we can apply POCS to $\mathscr{A} - \epsilon \mathscr{E}$ with a very small $\epsilon$.
If the iteration converges and both $\{ \mathscr{A}^{(t)} \}$ and $\{ \mathscr{B}^{(t)} \}$ converge to the same tensor, then we can conclude that $\mathscr{A}$ is M-PD, i.e., the strong ellipticity holds.

\section{Further sufficient conditions in particular cases}\label{sec_special}

One may conjecture that every M-PSD tensor can be modified into an S-PSD by preserving the summations of the corresponding entries.
Unfortunately, this can be disproved by the following counter example given by Choi and Lam \cite{Choi75,ChoiLam77}:
\begin{equation}\label{eq_choi}
\begin{split}
  \mathscr{A}_\gamma {\bf x}^2 {\bf y}^2 = &
  \ x_1^2 y_1^2 + x_2^2 y_2^2 + x_3^2 y_3^2 \\
  & - 2(x_1 x_2 y_1 y_2 + x_2 x_3 y_2 y_3 + x_3 x_1 y_3 y_1) \\
  & + \gamma (x_1^2 y_2^2 + x_2^2 y_3^2 + x_3^2 y_1^2),
\end{split}
\end{equation}
where $\gamma \geq 1$.
When we apply POCS to tensor $\mathscr{A}_\gamma$, the limitations of $\{ \mathscr{A}^{(t)} \}$ and $\{ \mathscr{B}^{(t)} \}$ differ from each other, which implies that the intersection of $\mathbb{T}_{\mathscr{A}_\gamma}$ and $\mathbb{S}$ is empty.
Furthermore, Choi and Lam \cite{ChoiLam77} also proved that $\{ \lambda \mathscr{A}_1:\, \lambda \geq 0 \}$ is an extremal ray of the M-PSD cone.
Therefore, we desire to extend our sufficient condition of the SE-condition to include more situations.

Given a fourth-order tensor $\mathscr{A} \in \mathbb{E}$, we can always calculate the eigendecomposition of its unfolding matrix
$$
{\bf A} = \sum_{s=1}^r \alpha_s {\bf u}_s {\bf u}_s^\top,
$$
thus matrix $\mathscr{A} {\bf y}^2$ can be correspondingly decomposed into
\begin{equation}\label{eq_grad}
  \mathscr{A} {\bf y}^2 = \sum_{s = 1}^r \alpha_s ({\bf U}_s {\bf y}) ({\bf U}_s {\bf y})^\top,
\end{equation}
where ${\bf u}_s$ is the vectorization of ${\bf U}_s$ ($s = 1,2,\dots,r$).
Note that the coefficients $\alpha_s$ are not necessarily positive, otherwise this is exactly the case discussed in the previous section.
The eigendecomposition of the unfolding matrix guarantees the existence of such decompositions as \eqref{eq_grad}.
Actually, the orthogonality of the vectorizations of ${\bf U}_s$ is not required in the following discussion, and the number of terms, i.e., $r$, may also be larger than nine.
Without loss of generality, we assume that $\alpha_1,\dots,\alpha_q > 0$ and $\alpha_{q+1},\dots,\alpha_r < 0$.

{\bf Case 1:} $q=3$ and ${\bf U}_s$ ($s=1,2,3$) are rank-one, i.e.,
$$
{\bf U}_s = {\bf v}_s {\bf w}_s^\top, \quad s = 1,2,3,
$$
and furthermore ${\bf V} = \big[ {\bf v}_1, {\bf v}_2, {\bf v}_3 \big]$ and ${\bf W} = \big[ {\bf w}_1, {\bf w}_2, {\bf w}_3 \big]$ are nonsingular.
Then the summation of the first three PSD terms is written into
$$
\sum_{s = 1}^3 \alpha_s ({\bf U}_s {\bf y}) ({\bf U}_s {\bf y})^\top =
{\bf V}
\begin{bmatrix}
  \alpha_1 ({\bf w}_1^\top {\bf y})^2 \\
  & \alpha_2 ({\bf w}_2^\top {\bf y})^2 \\
  & & \alpha_3 ({\bf w}_3^\top {\bf y})^2
\end{bmatrix}
{\bf V}^\top.
$$
Denote $\widehat{\bf U}_s := {\bf V}^{-1} {\bf U}_s$ and $\widehat{\bf U}_s = [\widehat{\bf u}_{s1},\widehat{\bf u}_{s2},\widehat{\bf u}_{s3}]^\top$ ($s = 4,5,\dots,r$).
If $\mathscr{A} {\bf y}^2$ is PSD for any nonzero ${\bf y} \in \mathbb{R}^3$, then $\sum_{s = 1}^4 \alpha_s ({\bf U}_s {\bf y}) ({\bf U}_s {\bf y})^\top$ is also PSD since $\alpha_5,\dots,\alpha_r < 0$.
That is,
$$
\begin{bmatrix}
  \alpha_1 ({\bf w}_1^\top {\bf y})^2 \\
  & \alpha_2 ({\bf w}_2^\top {\bf y})^2 \\
  & & \alpha_3 ({\bf w}_3^\top {\bf y})^2
\end{bmatrix}
+
\alpha_4
\begin{bmatrix}
  \widehat{\bf u}_{41}^\top {\bf y} \\ \widehat{\bf u}_{42}^\top {\bf y} \\ \widehat{\bf u}_{43}^\top {\bf y}
\end{bmatrix}
\begin{bmatrix}
  \widehat{\bf u}_{41}^\top {\bf y} & \widehat{\bf u}_{42}^\top {\bf y} & \widehat{\bf u}_{43}^\top {\bf y}
\end{bmatrix}
\succeq 0.
$$
When ${\bf y}$ is selected such that ${\bf w}_l^\top {\bf y} = 0$, it must hold that $\widehat{\bf u}_{sl}^\top {\bf y} = 0$ to guarantee the positive definiteness, which implies that $\widehat{\bf u}_{sl} = \sigma_{sl} {\bf w}_l$ ($l=1,2,3$).
Similarly, we can prove that in this case if $\mathscr{A}$ is M-PSD then
\begin{equation}\label{eq_case1_cond1}
  {\bf U}_s  = {\bf V} {\bm \Sigma}_s {\bf W}^\top, \quad s = 4,5,\dots,r,
\end{equation}
where ${\bm \Sigma}_s = {\rm diag}(\sigma_{s1},\sigma_{s2},\sigma_{s3})$.
Denote ${\bf D} = {\rm diag}({\bf w}_1^\top {\bf y},{\bf w}_2^\top {\bf y},{\bf w}_3^\top {\bf y})$ and ${\bm \sigma}_s = [\sigma_{s1},\sigma_{s2},\sigma_{s3}]^\top$ ($s = 4,5,\dots,r$).
Then $\mathscr{A} {\bf y}^2 = {\bf D} {\bf C} {\bf D}^\top$, where
\begin{equation}
{\bf C} =
{\rm diag}(\alpha_1, \alpha_2, \alpha_3)
+ \sum_{s=4}^r \alpha_s {\bm \sigma}_s {\bm \sigma}_s^\top.
\end{equation}
Therefore, $\mathscr{A}$ in this case is M-PSD if and only if ${\bf U}_s$ ($s = 4,5,\dots,r$) are in the form \eqref{eq_case1_cond1} and the matrix ${\bf C}$ is PSD.

\begin{theorem}
  Let $\mathscr{A} \in \mathbb{E}$ be given by \eqref{eq_grad} with (i) $\alpha_1,\alpha_2,\alpha_3 > 0$, $\alpha_4,\dots,\alpha_r < 0$, (ii) ${\bf U}_s = {\bf v}_s {\bf w}_s^\top$ ($s=1,2,3$), and (iii) ${\bf V} = \big[ {\bf v}_1, {\bf v}_2, {\bf v}_3 \big]$, ${\bf W} = \big[ {\bf w}_1, {\bf w}_2, {\bf w}_3 \big]$ are nonsingular. Then $\mathscr{A}$ is M-PSD if and only if (i) ${\bf U}_s  = {\bf V} {\bm \Sigma}_s {\bf W}^\top$ with ${\bm \Sigma}_s = {\rm diag}(\sigma_{s1},\sigma_{s2},\sigma_{s3})$ ($s = 4,\dots,r$), and (ii) the matrix ${\rm diag}(\alpha_1, \alpha_2, \alpha_3) + \sum_{s=4}^r \alpha_s {\bm \sigma}_s {\bm \sigma}_s^\top$ is PSD.
\end{theorem}

Nevertheless, we can show that the tensors of this type also satisfy the sufficient condition established in Section \ref{sec_convex}.
We can calculate the eigendecomposition of ${\bf C} = \sum_{s=1}^{\widetilde{r}} \widetilde{\alpha}_s \widetilde{\bm \sigma}_s \widetilde{\bm \sigma}_s^\top$.
Then all the coefficients $\widetilde{\alpha}_s$ are positive since ${\bf C}$ is PSD.
Denote $\widetilde{\bf U}_s := {\bf V} \widetilde{\bm \Sigma}_s {\bf W}^\top$, where $\widetilde{\bf \Sigma}_s = {\rm diag}(\widetilde{\sigma}_{s1},\widetilde{\sigma}_{s2},\widetilde{\sigma}_{s3})$ ($s = 1,2,\dots,\widetilde{r}$).
We can easily verify that $\mathscr{A} {\bf y}^2 = \sum_{s = 1}^{\widetilde{r}} \widetilde{\alpha}_s (\widetilde{\bf U}_s {\bf y}) (\widetilde{\bf U}_s {\bf y})^\top$ with $\widetilde{\alpha}_s > 0$ ($s = 1,2,\dots,\widetilde{r}$).

{\bf Case 2}: $r=7$, $q=6$, and ${\bf U}_s$ ($s = 1,2,\dots,6$) are rank-one with
$$
{\bf U}_s =
\left\{
\begin{array}{ll}
  {\bf v}_s {\bf w}_s^\top, & s = 1,2,3, \\
  {\bf v}_{s-3} {\bf w}_s^\top, & s = 4,5,6,
\end{array}
\right.
$$
and ${\bf V} = [{\bf v}_1, {\bf v}_2, {\bf v}_3]$, ${\bf W} = [{\bf w}_1, {\bf w}_2, {\bf w}_3]$, $\widetilde{\bf W} = [{\bf w}_4, {\bf w}_5, {\bf w}_6]$ are nonsingular.
It is reasonable to further assume that ${\bf w}_s$ and ${\bf w}_{s+3}$ are linearly independent ($s=1,2,3$).
Otherwise we can simply collect them into one term.
Similarly to our analysis of Case 1, we can conclude that $\mathscr{A}$ in Case 2 is M-PSD if and only if
$$
{\bf U}_7 = {\bf V} {\bm \Sigma}_7 {\bf W}^\top + {\bf V} \widetilde{\bm \Sigma}_7 \widetilde{\bf W}^\top,
$$
where ${\bm \Sigma}_7 = {\rm diag}(\sigma_{1},\sigma_{2},\sigma_{3})$ and $\widetilde{\bm \Sigma}_7 = {\rm diag}(\sigma_{4},\sigma_{5},\sigma_{6})$, and
\begin{multline}\label{eq_case2_cond1}
\begin{bmatrix}
  \alpha_1 ({\bf w}_1^\top {\bf y})^2 + \alpha_4 ({\bf w}_4^\top {\bf y})^2 \\
  & \alpha_2 ({\bf w}_2^\top {\bf y})^2 + \alpha_5 ({\bf w}_5^\top {\bf y})^2 \\
  & & \alpha_3 ({\bf w}_3^\top {\bf y})^2 + \alpha_6 ({\bf w}_6^\top {\bf y})^2
\end{bmatrix} \\
+
\alpha_7
\begin{bmatrix}
  (\sigma_{1} {\bf w}_{1} + \sigma_{4} {\bf w}_4)^\top {\bf y} \\ (\sigma_{2} {\bf w}_{2} + \sigma_{5} {\bf w}_5)^\top {\bf y} \\ (\sigma_{3} {\bf w}_{3} + \sigma_{6} {\bf w}_6)^\top {\bf y}
\end{bmatrix}
\cdot
\begin{bmatrix}
  (\sigma_{1} {\bf w}_{1} + \sigma_{4} {\bf w}_4)^\top {\bf y} \\ (\sigma_{2} {\bf w}_{2} + \sigma_{5} {\bf w}_5)^\top {\bf y} \\ (\sigma_{3} {\bf w}_{3} + \sigma_{6} {\bf w}_6)^\top {\bf y}
\end{bmatrix}^\top
\succeq 0.
\end{multline}
Denote
\begin{equation}\label{eq_case2_fun}
  \eta({\bf y}) := \sum_{s=1}^3 \frac{(\sigma_s {\bf w}_s^\top {\bf y} + \sigma_{s+3} {\bf w}_{s+3}^\top {\bf y})^2}{\alpha_s ({\bf w}_s^\top {\bf y})^2 + \alpha_{s+3} ({\bf w}_{s+3}^\top {\bf y})^2}.
\end{equation}
Then the matrix in \eqref{eq_case2_cond1} is PSD for all nonzero ${\bf y}$ if and only if
\begin{equation}\label{eq_case2_cond2}
  \sup \bigg\{ \eta({\bf y}):\, {\bf y} \notin \bigcup_{s=1,2,3} \big( {\bf w}_s^\bot \cap {\bf w}_{s+3}^\bot \big) \bigg\} \leq \frac{1}{-\alpha_7},
\end{equation}
where ${\bf w}_s^\bot$ denotes the orthogonal complement subspace to ${\rm span}({\bf w}_s)$.
Furthermore, the elasticity tensor in both Case 1 and Case 2 cannot be M-PD, which can be seen by taking a vector ${\bf y}$ in $\bigcup_{s=1,2,3} \big( {\bf w}_s^\bot \cap {\bf w}_{s+3}^\bot)$ and thus $\mathscr{A} {\bf y}^2$ has a nonempty null space.

\begin{theorem}
  Let $\mathscr{A} \in \mathbb{E}$ be given by \eqref{eq_grad} with (i) $\alpha_1,\dots,\alpha_6 > 0$, $\alpha_7 < 0$, (ii) ${\bf U}_s = {\bf v}_s {\bf w}_s^\top$ ($s=1,2,\dots,6$), ${\bf v}_s = {\bf v}_{s+3}$ ($s=1,2,3$), (iii) ${\bf V} = \big[ {\bf v}_1, {\bf v}_2, {\bf v}_3 \big]$, ${\bf W} = \big[ {\bf w}_1, {\bf w}_2, {\bf w}_3 \big]$, $\widetilde{\bf W} = \big[ {\bf w}_4, {\bf w}_5, {\bf w}_6 \big]$ are nonsingular, and (iv) ${\bf w}_s$ and ${\bf w}_{s+3}$ are linearly independent ($s=1,2,3$). Then $\mathscr{A}$ is M-PSD if and only if (i) ${\bf U}_7 = {\bf V} {\bm \Sigma}_7 {\bf W}^\top + {\bf V} \widetilde{\bm \Sigma}_7 \widetilde{\bf W}^\top$ with ${\bm \Sigma}_7 = {\rm diag}(\sigma_{1},\sigma_{2},\sigma_{3})$ and $\widetilde{\bm \Sigma}_7 = {\rm diag}(\sigma_{4},\sigma_{5},\sigma_{6})$, and (ii) $\sup \Big\{ \eta({\bf y}):\, {\bf y} \notin \bigcup_{s=1,2,3} \big( {\bf w}_s^\bot \cap {\bf w}_{s+3}^\bot \big) \Big\} \leq \frac{1}{-\alpha_7}$, where $\eta(\cdot)$ is defined by \eqref{eq_case2_fun}.
\end{theorem}

Generally speaking, the supreme in the left-hand side of \eqref{eq_case2_cond2} is not easy to obtain. In some particular case, however, we can calculate the exact supreme of $\eta({\bf y})$ in the domain of definition. For instance, recall the counter example given in \eqref{eq_choi}, where the matrix $\mathscr{A}_1 {\bf y}^2$ can be decomposed into
\[
\begin{split}
\mathscr{A}_1 {\bf y}^2
&=
\begin{bmatrix}
  y_1^2 + y_2^2 & -y_1 y_2 & -y_3 y_1 \\
  -y_1 y_2 & y_2^2 + y_3^2 & -y_2 y_3 \\
  -y_3 y_1 & -y_2 y_3 & y_3^2 + y_1^2
\end{bmatrix} \\
&=
\begin{bmatrix}
  2 y_1^2 + y_2^2 \\
  & 2 y_2^2 + y_3^2 \\
  & & 2 y_3^2 + y_1^2
\end{bmatrix}
-
\begin{bmatrix}
  y_1 \\ y_2 \\ y_3
\end{bmatrix}
\cdot
\begin{bmatrix}
  y_1 & y_2 & y_3
\end{bmatrix}.
\end{split}
\]
In this example, the coefficient $\alpha_7 = -1$ and the function
$$
\eta({\bf y}) = \frac{y_1^2}{2 y_1^2 + y_2^2} + \frac{y_2^2}{2 y_2^2 + y_3^2} + \frac{y_3^2}{2 y_3^2 + y_1^2}.
$$
We can verify that $\eta({\bf y}) \leq 1$ and the equality holds if $y_1 = y_2 = y_3 \neq 0$, which implies $\sup \big\{ \eta({\bf y}):\, (y_1,y_2)\neq(0,0), (y_2,y_3)\neq(0,0), (y_3,y_1)\neq(0,0) \big\} = 1 = -1/\alpha_7$.
Therefore $\mathscr{A}_1$ is M-PSD.

{\bf Case 3}: $r=10$, $q=9$, and ${\bf U}_s$ ($s = 1,2,\dots,9$) are rank-one with
$$
{\bf U}_s =
\left\{
\begin{array}{ll}
  {\bf v}_s {\bf w}_s^\top, & s = 1,2,3, \\
  {\bf v}_{s-3} {\bf w}_s^\top, & s = 4,5,6, \\
  {\bf v}_{s-6} {\bf w}_s^\top, & s = 7,8,9,
\end{array}
\right.
$$
and ${\bf V} = [{\bf v}_1, {\bf v}_2, {\bf v}_3]$, ${\bf W} = [{\bf w}_1, {\bf w}_2, {\bf w}_3]$, $\widetilde{\bf W} = [{\bf w}_4, {\bf w}_5, {\bf w}_6]$, $\widehat{\bf W} = [{\bf w}_7, {\bf w}_8, {\bf w}_9]$ are nonsingular.
It is also reasonable to further assume that $\{ {\bf w}_s, {\bf w}_{s+3}, {\bf w}_{s+6} \}$ are linearly independent ($s=1,2,3$).
Similarly to Cases 1 and 2, we can conclude that $\mathscr{A}$ in Case 3 is M-PSD if and only if
$$
{\bf U}_{10} = {\bf V} {\bm \Sigma}_{10} {\bf W}^\top + {\bf V} \widetilde{\bm \Sigma}_{10} \widetilde{\bf W}^\top + {\bf V} \widehat{\bm \Sigma}_{10} \widehat{\bf W}^\top,
$$
where ${\bm \Sigma}_{10} = {\rm diag}(\sigma_{1},\sigma_{2},\sigma_{3})$, $\widetilde{\bm \Sigma}_{10} = {\rm diag}(\sigma_{4},\sigma_{5},\sigma_{6})$, $\widehat{\bm \Sigma}_{10} = {\rm diag}(\sigma_{7},\sigma_{8},\sigma_{9})$, and
\begin{multline}\label{eq_case3_cond1}
\begin{bmatrix}
  \sum\limits_{h=0,1,2} \alpha_{1+3h} ({\bf w}_{1+3h}^\top {\bf y})^2 \\
  & \sum\limits_{h=0,1,2} \alpha_{2+3h} ({\bf w}_{2+3h}^\top {\bf y})^2 \\
  & & \sum\limits_{h=0,1,2} \alpha_{3+3h} ({\bf w}_{3+3h}^\top {\bf y})^2
\end{bmatrix} \\
+
\alpha_{10}
\begin{bmatrix}
  (\sigma_{1} {\bf w}_1 + \sigma_{4} {\bf w}_4 + \sigma_{7} {\bf w}_7)^\top {\bf y} \\ (\sigma_{2} {\bf w}_2 + \sigma_{5} {\bf w}_5 + \sigma_{8} {\bf w}_8)^\top {\bf y} \\ (\sigma_{3} {\bf w}_3 + \sigma_{6} {\bf w}_6 + \sigma_{9} {\bf w}_9)^\top {\bf y}
\end{bmatrix}
\cdot
\begin{bmatrix}
  (\sigma_{1} {\bf w}_1 + \sigma_{4} {\bf w}_4 + \sigma_{7} {\bf w}_7)^\top {\bf y} \\ (\sigma_{2} {\bf w}_2 + \sigma_{5} {\bf w}_5 + \sigma_{8} {\bf w}_8)^\top {\bf y} \\ (\sigma_{3} {\bf w}_3 + \sigma_{6} {\bf w}_6 + \sigma_{9} {\bf w}_9)^\top {\bf y}
\end{bmatrix}^\top
\succeq 0.
\end{multline}
Denote
\begin{equation}\label{eq_case3_fun}
  \eta({\bf y}) := \sum_{s=1}^3 \frac{(\sigma_s {\bf w}_s^\top {\bf y} + \sigma_{s+3} {\bf w}_{s+3}^\top {\bf y} + \sigma_{s+6} {\bf w}_{s+6}^\top {\bf y})^2}{\alpha_s ({\bf w}_s^\top {\bf y})^2 + \alpha_{s+3} ({\bf w}_{s+3}^\top {\bf y})^2 + \alpha_{s+6} ({\bf w}_{s+6}^\top {\bf y})^2}.
\end{equation}
According to the assumptions, $\eta({\bf y})$ is well-defined for any nonzero ${\bf y}$.
Thus the matrix in \eqref{eq_case3_cond1} is PSD for all nonzero ${\bf y}$ if and only if
$$
\max \big\{ \eta({\bf y}):\, {\bf y}^\top {\bf y} = 1 \big\} \leq \frac{1}{-\alpha_{10}}.
$$
Furthermore, Case 3 is of interest since ${\bf V} = [{\bf e}_1,{\bf e}_2,{\bf e}_3]$, ${\bf W} = [{\bf e}_1,{\bf e}_2,{\bf e}_3]$, $\widetilde{\bf W} = [{\bf e}_2,{\bf e}_3,{\bf e}_1]$, $\widehat{\bf W} = [{\bf e}_3,{\bf e}_1,{\bf e}_2]$ are exactly the case for isotropic and some particular anisotropic linearly elastic materials \cite{ZubovRudev16}.

\begin{theorem}
  Let $\mathscr{A} \in \mathbb{E}$ be given by \eqref{eq_grad} with (i) $\alpha_1,\dots,\alpha_9 > 0$, $\alpha_{10} < 0$, (ii) ${\bf U}_s = {\bf v}_s {\bf w}_s^\top$ ($s=1,2,\dots,9$), ${\bf v}_s = {\bf v}_{s+3} = {\bf v}_{s+6}$ ($s=1,2,3$), (iii) ${\bf V} = \big[ {\bf v}_1, {\bf v}_2, {\bf v}_3 \big]$, ${\bf W} = \big[ {\bf w}_1, {\bf w}_2, {\bf w}_3 \big]$, $\widetilde{\bf W} = \big[ {\bf w}_4, {\bf w}_5, {\bf w}_6 \big]$, $\widehat{\bf W} = \big[ {\bf w}_7, {\bf w}_8, {\bf w}_9 \big]$ are nonsingular, and (iv) ${\bf w}_s, {\bf w}_{s+3}, {\bf w}_{s+6}$ are linearly independent ($s=1,2,3$). Then $\mathscr{A}$ is M-PSD if and only if (i) ${\bf U}_{10} = {\bf V} {\bm \Sigma}_{10} {\bf W}^\top + {\bf V} \widetilde{\bm \Sigma}_{10} \widetilde{\bf W}^\top + {\bf V} \widehat{\bm \Sigma}_{10} \widehat{\bf W}^\top$ with ${\bm \Sigma}_{10} = {\rm diag}(\sigma_{1},\sigma_{2},\sigma_{3})$, $\widetilde{\bm \Sigma}_{10} = {\rm diag}(\sigma_{4},\sigma_{5},\sigma_{6})$, $\widehat{\bm \Sigma}_{10} = {\rm diag}(\sigma_{7},\sigma_{8},\sigma_{9})$, and (ii) $\max \big\{ \eta({\bf y}):\, {\bf y}^\top {\bf y} = 1 \big\} \leq \frac{1}{-\alpha_{10}}$, where $\eta(\cdot)$ is defined by \eqref{eq_case3_fun}. Furthermore, $\mathscr{A}$ is M-PD when the strict inequality holds in the last condition.
\end{theorem}

\section{Conclusions}\label{sec_conclusion}

We have established several sufficient conditions for the strong ellipticity (M-positive definiteness) of general elasticity tensors.
S-positive definiteness has been already known as an easily checkable sufficient conditions for M-positive definiteness.
However, the range of S-positive (semi)definiteness is too narrow.
Thus our first sufficient condition extends the coverage of S-PSD tensors, which states that $\mathscr{A}$ is M-PSD or M-PD if it can be modified into an S-PSD or S-PD tensor $\mathscr{B}$ respectively by preserving $b_{ijkl} = b_{jilk}$ and $b_{ijkl} + b_{jikl} = a_{ijkl} + a_{jikl}$.
To check whether a tensor satisfies this condition, we employ an alternating projection method called POCS and verify its convergence.

Next, we have considered three particular cases and provided the necessary and sufficient conditions of the strong ellipticity.
Actually, we can understand these three cases in the following way.
Any tensor $\mathscr{A} \in \mathbb{E}$ can be factorized into two parts $\mathscr{A} = \mathscr{A}_1 - \mathscr{A}_2$, where $\mathscr{A}_1$ and $\mathscr{A}_2$ satisfy that $\mathbb{T}_{\mathscr{A}_1} \cap \mathbb{S} \neq \emptyset$ and $\mathbb{T}_{\mathscr{A}_2} \cap \mathbb{S} \neq \emptyset$.
We further assume that $\mathscr{A}_1 = \sum_{s=1}^9 \alpha_s {\bf v}_s \circ {\bf v}_s \circ {\bf w}_s \circ {\bf w}_s$ ($\alpha_s \geq 0$), where $\circ$ stands for the outer product of tensors.
It is interesting to note that such $\mathscr{A}_1$ is in the dual cone of the M-PSD cone.
We have established necessary and sufficient conditions for the strong ellipticity for the cases in which ${\bf v}_s$ and ${\bf w}_s$ satisfy some dependence conditions and the rank of $\mathscr{A}_2 {\bf y}^2$ is no greater than one for all ${\bf y} \in \mathbb{R}^3$. 

Our further investigation will be threefold.
The first one is to find out whether there are other M-PSD or M-PD elasticity tensors beyond our sufficient conditions.
Alternatively, we desire to know whether every M-PSD elasticity tensor can be represented by a convex combination of several tensors satisfying the conditions presented in this paper.
The next part is the characterization of the M-PSD cone. We already know that $\{ \alpha \mathscr{A}:\, \alpha \geq 0 \}$ is an extremal ray of the M-PSD cone if there is a tensor in $\mathbb{T}_{\mathscr{A}} \cap \mathbb{S}$ whose unfolding matrix is rank-one.
We would like to find more extremal rays of this convex cone.
Of course, if all the extremal rays were discovered, then the M-PSD cone could be characterized thoroughly.
Our third target in the future is to generalize these results to higher order elasticity tensors.


\end{document}